\documentclass[12pt]{article}
\usepackage{amssymb,amsmath}
\usepackage{latexsym, amsfonts, amsmath}%, wasysym}

\newcommand{\qed}{\hskip 5mm \rule{2.5mm}{2.5mm}\vskip 10pt}
\newcommand{\R}{{\mathbb R}}

\newcommand{\N}{{\mathbb N}}

\newcommand{\proof}{\noindent{\em Proof:\ }}

\newcommand{\var}{\mbox{\rm var}}

\begin{document}
\newtheorem{thm}{Theorem}[section]
\newtheorem{defs}[thm]{Definition}
\newtheorem{lem}[thm]{Lemma}
\newtheorem{note}[thm]{Note}
\newtheorem{rem}[thm]{Remark}
\newtheorem{cor}[thm]{Corollary}
\newtheorem{prop}[thm]{Proposition}
\renewcommand{\theequation}{\arabic{section}.\arabic{equation}}
\newcommand{\newsection}[1]{\setcounter{equation}{0} \section{#1}}
%%%%%%%%%%%% title %%%%%%%%%%%%%%%%%%%%%%%%%%%%%%%%
\title{Bernoulli Processes in Riesz spaces
   }
%%%%%%%%%%%%%%%%%%%%%%%%%%%%%%%%%%%%%%%%%%%%%%%%
\author{
 Wen-Chi Kuo\footnote{Research conducted while on an NRF post doctoral fellowship}\\
School of Computational and Applied Mathematics\\
 University of the Witwatersrand\\
 Private Bag 3, P O WITS 2050, South Africa \\ \\
 Jessica Vardy\footnote{Funded in part by the FRC }\  \&
 Bruce A. Watson\footnote{Funded in part by NRF grant  IFR2011032400120  and the Centre for Applicable 
Analysis and Number Theory} \\ 
 School of Mathematics\\
 University of the Witwatersrand\\
 Private Bag 3, P O WITS 2050, South Africa }
\maketitle
%%%%%%%%%%%%
\abstract{The action and averaging properties of conditional expectation operators
 are studied in the, measure-free, Riesz space,
 setting of Kuo, Labuschagne and Watson
[{ Conditional expectations on Riesz spaces},
          {\em J. Math. Anal. Appl.}, {\bf 303} (2005), 509-521]
 but on the abstract
 $L^2$ space, ${\cal L}^2(T)$ introduced by 
 Labuschagne and Watson 
         [{ Discrete Stochastic Integration in Riesz Spaces},
          {\em Positivity}, {\bf 14}, (2010), 859 - 575].
 In this setting it is shown that conditional expectation operators leave  ${\cal L}^2(T)$ invariant
 and  the Bienaym\'e equality and Tchebichev inequality are proved. 
 From this foundation Bernoulli processes are considered.
  Bernoulli's strong law of large numbers  and Poisson's theorem are formulated and proved.
}
%%%%%%%%%%%
\parindent=0cm
\parskip=0.5cm
%%%%%%%%%%%%%%%%%%%%%%%%%%%%%%%%%%
\newsection{Introduction}

Various authors have considered generalizations of stochastic processes to 
vector lattices / Riesz spaces, with a variety of assumptions being made on the 
processes being considered.  Most of this work has focussed on martingale theory, see, for 
example, \cite{demarr}, \cite{grobler-cts}, \cite{klw-indag}, \cite{lw}, \cite{stoica-4} 
and \cite{troitski}. The abstract properties of conditional expectation operators have also been
explored in various settings, see \cite{dodds}, \cite{klw-exp}, \cite{LdP} and \cite{schaefer}
and \cite{watson}.  However, the more elementary processes such as Markov processes, 
see \cite{vw-1}, Bernoulli processes and Poisson processes, which just rely only on  the concepts
of a conditional expecation operator and independence, have received  little attention.
As these processes have less accessible structure, their
study relies more heavily on properties of the underlying Riesz space, 
the representation of the conditional expectation operators and multiplication operations
in Riesz spaces.  If a Riesz space has a weak order unit, then
the order ideal generated by a weak order unit is order dense in the space.  However the order
ideal generated by the weak order unit is an $f$-algebra, see \cite{BBT}, \cite{BvR} and 
\cite{zaanen}, giving a multiplicative structure on a dense subspace.  Much of the
work in this paper lies on a Riesz space vector analogue of $L^2$ and the action of 
conditional expectations in this space, see Theorem \ref{lemmastar}, and their averaging
property, see  Lemma \ref{lem-lw}.
In particular the Bienaym\'e equality,
Theorem \ref{bienayme}, will be posed in 
this setting.  The Bienaym\'e equality enables us to give a Riesz space
analogue of Bernoulli's law of large numbers, Theorem \ref{bernoulli}.
One important property of the ideal generated by the weak order unit is that it posesses a 
functional calculus which enables one to lift continuous real valued functions on $[0,1]$ to 
the Riesz space, see \cite{grobler2} and \cite{zaanen}.  This is critical for Poisson's theorem,
Theorem \ref{poisson}.  We refer the reader to \cite{loeve} for the classical version of  
the Bienaym\'e equality, the Bernoulli law of large numbers and Poisson's theorem.

%%%%%%%%%%%%%%%%%%%%%%%%%%%
\newsection{Riesz Space Preliminaries}
 
We refer the reader to \cite{A-B} and \cite{zaanen} for general Riesz space theory.
The definitions and preliminaries presented here are specific to Riesz spaces with a weak order
unit and a conditional expectation operator.

The notion of a conditional expectation operator in a Dedekind complete Riesz space, $E$, with 
weak order unit was introduced in \cite{klw-indag} as a
 positive order continuous projection $T:E\to E,$ with
 range ${\cal R}(T)$ a Dedekind complete Riesz subspace of $E$, and having 
 $Te$ a weak order unit of $E$ for each weak order unit $e$ of $E$. 
 Instead of requiring $Te$ to be a weak order unit of $E$ for each weak order unit $e$ of $E$
 one can equivalently
 impose that there is a weak order unit in $E$ which is invariant under $T$.
 Averaging properties of conditional expectation operators and various other structural
 aspects were considered in \cite{klw-exp}.  In particular if 
  $B$ is the band in $E$ generated by $0\leq g\in {\cal R}(T)$
 and $P$ is the band projection onto $B$, it was shown that
 $Tf\in B$, for each $f\in B$, $Pf, (I-P)f \in{\cal R}(T)$
 for each $f\in {\cal R}(T)$, where $I$ denotes the identity map, and
 $Tf\in B^d$, for each $f\in B^d$.  A consequence of these relations and Freudenthal's
 theorem, \cite{zaanen}, is that 
 if $B$ is the band in $E$ 
  generated by $0\le g\in {\cal R}(T)$, with associated band projection $P$, then $TP=PT$,
 see \cite{klw-exp} for details.

To access the averaging properties of conditional expectation operators a multiplicative structure is
needed.  In the Riesz space setting the most natural multiplicative structure is that of
an $f$-algebra. This gives a multiplicative structure that is compatible with 
the order and additive structures on the space.  
The ideal, $E^e$, of $E$ generated by $e$, where $e$ is a weak order
unit of $E$ and $E$ is Dedekind complete, has a natural $f$-algebra structure.
This is constructed by setting
$(Pe)\cdot (Qe)=PQe=(Qe)\cdot (Pe)$ for band projections $P$ and $Q$,
and extending to $E^e$ by use of Freudenthal's Theorem.
In fact this process extends the multiplicative structure to the universal completion, $E^u$, of $E$.
This multiplication is associative, 
distributive and is positive in the sense that if $x,y\in E^+$ then $xy\ge 0$.
Here $e$ is the multiplicative unit.
For more information about $f$-algebras  see \cite{BBT, BvR, dodds, grobler, klw-exp, zaanen}.
If $T$ is a conditional expectation operator on the Dedekind complete Riesz space $E$ with weak 
order unit $e  = Te$, then restricting our attention to the $f$-algebra $E^e$,
$T$ is an averaging operator on $E^e$ if $T(fg)=fTg$ for $f,g\in E^e$ and $f\in R(T)$, see 
\cite{dodds, grobler, klw-exp}.  More will said about averaging operators in Section 3.

In a Dedekind complete Riesz space, $E$, with weak order unit and $T$
 a strictly positive conditional expectation on $E$.
 We say that the space is $T$-universally complete 
 if for each increasing net $(f_\alpha)$ in $E_+$
 with $(Tf_\alpha)$ order bounded in the universal completion $E^u$, 
 we have that $(f_\alpha)$ is order convergent.  If this is not the case, then both the space
 and conditional expectation operator can be extended so that the extended space is 
 $T$-universally complete with respect to the extended $T$, see \cite{klw-exp}.
 The extended space is also know as the natural domain of $T$, denoted ${\rm dom}(T)$ 
  of ${\cal L}^1(T)$, see \cite{dodds, grobler}.

 Let $E$ be a Dedekind complete Riesz space with conditional expectation $T$ and weak order unit $e=Te$.
 If $P$ and $Q$ are band projections on $E$, we say that $P$ and $Q$ are $T$-conditionally independent with respect to $T$ if
 \begin{eqnarray}
  TPTQe=TPQe=TQTPe.\label{indep-e}
 \end{eqnarray}
 We say that two Riesz subspaces $E_1$ and $E_2$ of $E$ are $T$-conditionally independent with respect to $T$ if
 all band projections $P_i, i=1,2,$ in $E$ with  $P_ie\in E_i, i=1,2,$ are $T$-conditionally independent with respect to $T$.
Equivalently (\ref{indep-e}) can be replaced with
 \begin{eqnarray}
  TPTQw=TPQw=TQTPw\quad\mbox{for all}\quad w\in \mathcal{R}(T).\label{indep-w}
 \end{eqnarray}
It should be noted that 
 $T$-conditional independence of the band projections $P$ and $Q$ is equivalent to
$T$-conditional independence of the closed Riesz subspaces $\left< Pe, \mathcal{R}(T)\right>$ and $\left< Qe, \mathcal{R}(T)\right>$
generated by $Pe$ and $\mathcal{R}(T)$ and by $Qe$ and $\mathcal{R}(T)$ respectively.
From the Radon-Nikod\'ym-Douglas-And\^o type theorem was established in \cite{watson},
if $E$ is a $T$-universally complete,
a subset $F$ of $E$ is a closed Riesz subspace of $E$ with ${\cal R}(T)\subset F$ if and only if
 there is a unique conditional expectation $T_F$ on $E$ with
 ${\cal R}(T_F)=F$ and $TT_F=T=T_FT$. In this case $T_Ff$ for $f\in E^+$ is uniquely determined by the
 property that
 \begin{eqnarray}\label{R-N}
 TPf=TPT_Ff
 \end{eqnarray}
 for all band projections on $E$ with $Pe\in F$.
 As a consequence of this,
 two closed Riesz subspaces $E_1$ and $E_2$ with ${\cal R}(T)\subset E_1\cap E_2$ are
 $T$-conditionally independent,  if and only if 
 \begin{eqnarray}
    T_1T_2=T=T_2T_1,\label{indep-2014}
  \end{eqnarray}
 where $T_i$ is the conditional expectation commuting with $T$ and having range $E_i, i=1,2$.
 Here (\ref{indep-2014}) can be equivalently replaced by 
 \begin{eqnarray}
	T_if=Tf,\quad\mbox{for all}\quad f\in E_{3-i}, \quad i=1,2,\label{eq-indep-1}
 \end{eqnarray}
 see \cite{vw-1}.
The concept of $T$-conditional independence can be extended to a family, say
 $(E_\lambda)_{\lambda\in \Lambda}$, of closed Dedekind complete Riesz
 subspaces of $E$ with ${\cal R}(T)\subset E_\lambda$ for all $\lambda\in\Lambda$.
 We say that the family is $T$-conditionally independent
 if, for each pair of disjoint sets $\Lambda_1, \Lambda_2 \subset \Lambda$, we have that
 $E_{\Lambda_1}$ and $E_{\Lambda_2}$ are $T$-conditionally independent.
  Here $E_{\Lambda_j}: = \left< \bigcup_{\lambda\in\Lambda_j} E_\lambda \right>$.
  Finally, we say that a sequence $(f_n)$ in $E$ is $T$-conditionally independent if
 the family of closed Riesz subspaces $\left<\{f_n\}\cup \mathcal{R}(T)\right>,  n\in\N,$ 
 is $T$-conditionally independent.

%%%%%%%%%%%%%%%%%%%%%%%%%%%
\newsection{Conditional expectation operators in ${\cal L}^2(T)$}

In this work we assume that $E$ is $T$-universally complete and in this 
case we have ${\mathcal L}^1(T) = E$, see \cite{lw}.  As $E^u$, the universal completion of $E$,
is an $f$-algebra, multiplication of elements of $E$ is defined but does not necessarily result
in an element of $E$. This leads us, as in \cite{lw}, to define  
    \[{\mathcal L}^2(T):= \left\{ x\in {\mathcal L}^1(T) \vert x^2 \in {\mathcal L}^1(T)\right\}.\]
If $f,g\in{\mathcal L}^2(T)$ then in the
$f$-algebra $E^u$, $0\le (f\pm g)^2=f^2\pm 2fg+g^2$.  Thus  $\pm 2fg\le f^2+g^2$
and $2|fg|\le f^2+g^2\in{\cal L}^1(T)=E$. Hence $fg\in {\cal L}^1(T)=E$. As noted
in \cite{lw}, a consequence of this is 
that ${\cal L}^2(T)$ is a vector space.

The averaging property of conditional expectation operators only makes sense if it can be ensured 
that the all products involved remain in the space.  Theorem 4.3 of \cite{klw-exp} states that if 
$E$ is a Dedekind complete Riesz space with weak order unit, $T$ is a 
conditional expectation operator on $E$ and $E$ is also an $f$-algebra, then $T$ is an averaging 
operator, i.e. $T(fg)=gTf$ for all $f\in E, g\in{\cal R}(T)$. The averaging property is revisited in
\cite[Theorem 2.1]{lw} without proof. 
The variant of \cite[Theorem 2.1]{lw} drops the assumption that $E$ is an $f$-algebra, 
but imposes the additional conditions that $fg\in E$ and that $E$ is 
$T$-universally complete. A strengthened version of this is proved in Lemma \ref{lem-lw}.
This however does not address whether $Sf\in {\cal L}^2(T)$
for $f\in {\cal L}^2(T)$ and $S$ a conditional expectation operator on $E$ with 
$TS=T=ST$.  For this see Theorem \ref{lemmastar} below. 
As a consequence of Lemma \ref{lem-lw} and Theorem \ref{lemmastar},
we are able to conclude, see Theorem~\ref{cor-average}
below, that for such a conditional expectation operator, $S$, 
$S(fTg)=Tg\cdot Sf$ for all $f,g\in {\cal L}^2(T)$.

\begin{lem}\label{lem-lw}
Let $E$ be a Dedekind complete Riesz space with weak order unit, $e$, and $T$ is a
 conditional expectation operator on $E$ with $Te = e$.
 If $f,g,fg\in E$ with $g\in {\cal R}(T)$ then $g\cdot Tf\in E$ and $T(fg)=g\cdot Tf.$ 
\end{lem}

\proof
 \underline{\bf Case I: $f,g,fg \in E_+$ with $g\in  {\cal R}(T)$}
 Let $f_n=f\wedge ne$ and $g_n=g\wedge ne$.  Then $f_n\uparrow f$ and
 $g_n\uparrow g$. Here $f_n, g_n\in E^e_+$ with $g_n\in {\cal R}(T)$, so 
 \cite[Theorem 4.3]{klw-exp} can be applied to give
    $T(f_ng_m)=g_mT(f_n), m,n \in\N.$
 Thus 
 \begin{eqnarray}
  g_mT(f_n)=T(f_ng_m)\le T(fg),\quad m,n\in\N.
  \label{average-2014-1}
 \end{eqnarray} 
 Here $f_ng_m\uparrow fg$ in $E$, so from the order continuity of $T$,
 $T(f_ng_m)\uparrow  T(fg)$ in $E$.
 In the universal completion, $E^u$, of $E$, we have $g_mT(f_n)\uparrow gT(f)$,
 however, from (\ref{average-2014-1}), $g_mT(f_n)$ is bounded above by $T(fg)\in E$, so
 $g_mT(f_n)\uparrow gT(f)$ in $E$, giving $T(fg)=gT(f)$.

\underline{\bf Case II: $f,g,fg \in E$ with $g\in  {\cal R}(T)$}
From Case I,   $T(f^\pm g^\mp)=g^\mp T(f^\pm)$ and
  $T(f^\pm g^\pm)=g^\pm T(f^\pm)$, from which the result follows.
\qed

\begin{thm}\label{lemmastar}
Let $E$ be a $T$-universally complete Riesz space with weak order unit, $e$, where $T$ is a
strictly positive conditional expectation operator with
$Te = e$ and let $S$ be a conditional expectation operator on $E$ with $TS=T=ST$.
If $f \in {\mathcal L}^2(T)$ then $Sf \in {\mathcal L}^2(T)$.
\end{thm}

\proof
Let $f \in {\mathcal L}^2(T)$ and define $f_n = (ne\wedge \vert f\vert)\in E^e_+$, $n \in \N$.  
We note that  $E^e$ is an $f$-algebra.  
Here $Sf_n \in E^e_+$ and $f_n - Sf_n \in E^e$ and so $(f_n - Sf_n)^2 \in E^e_+$.  But 
\begin{align*}
(f_n - Sf_n)^2 &=f_n^2 - 2f_n\cdotp Sf_n + (Sf_n)^2.
\end{align*}
Thus,  
\begin{eqnarray}
  0\leq S(f_n - Sf_n)^2 = Sf_n^2 - 2S(f_n\cdotp Sf_n) +S (Sf_n)^2.\label{baw-2014-8}
\end{eqnarray}
As conditional expectation operators on $E^e$ are averaging operators, see \cite{klw-exp}, 
for each $g\in E^e$ we have 
\begin{align}
S(g\cdot Sf_n) &=Sf_n \cdotp Sg.\label{baw-2014-7}
\end{align}
Taking $g=f_n$ in (\ref{baw-2014-7}) gives
\begin{eqnarray}
 S(f_n\cdot Sf_n) =Sf_n \cdotp Sf_n=(Sf_n)^2,\label{baw-2104-9}
\end{eqnarray}
while taking $g=Sf_n$ gives
\begin{eqnarray}
 S(Sf_n\cdot Sf_n) =Sf_n \cdotp S(Sf_n)=(Sf_n)^2,\label{baw-2014-10}
\end{eqnarray}
 as $S$ is a projection.
Combining (\ref{baw-2104-9}) and (\ref{baw-2014-10}) with (\ref{baw-2014-8}) gives
\begin{eqnarray}
   Sf_n^2 \ge (Sf_n)^2.\label{baw-2014-11}
\end{eqnarray}
Since $f_n \uparrow \vert f\vert $ we have
 $f_n^2 \uparrow \vert f\vert^2$ and $Sf_n^2 \uparrow S\vert f\vert^2$, as $n\to\infty$,
from the order continuity of $S$.
Similarly $Sf_n\uparrow S|f|$ giving $(Sf_n)^2\uparrow (S|f|)^2$ as $n\to\infty$.
Hence taking $n\to\infty$ in (\ref{baw-2014-11}) yields
	\[(S\vert f \vert)^2 \leq Sf^2.\]
But $\vert Sf\vert \le S\vert f \vert$ so  
	\[(Sf)^2 \le Sf^2 \in E,\]
giving  $Sf  \in {\mathcal L}^2(T)$.
\qed

\begin{cor}\label{cor-average}
Let $E$ be a $T$-universally complete Riesz space with weak order unit, $e$, where $T$ is a
strictly positive conditional expectation operator with
$Te = e$. Let $S, J$ be conditional expectation operators on $E$ with $TS=T=ST$, $TJ=T=JT$
 and $JS=J=SJ$.
If $f, g \in {\mathcal L}^2(T)$ then $S(f\cdot Jg)=Jg\cdot S(f)$.
 \end{cor}

\proof
 As $f, g \in {\mathcal L}^2(T)$ from Theorem \ref{lemmastar} $Jg\in {\mathcal L}^2(T)$.
 Now $f, Jg \in {\mathcal L}^2(T)$ so $f, Jg, f\cdot Jg \in E$ so Lemma \ref{average-2014-1}
 gives $Jg\cdot Sf\in E$ and $S(f\cdot Jg)=Jg\cdot Sf$.
 \qed

We are now in a position to give the Tchebichev's inequality in ${\cal L}^2(T)$.

\begin{thm}[Tchebichev's Inequality]
 Let $E$ be a Dedekind complete Riesz space with conditional expectation $T$ and
 weak order unit $e=Te$.  Let $f\in {\cal L}^2(T), f\ge 0,$ and $\epsilon\in\R, \epsilon>0,$
 then
 $$TP_{(f-\epsilon e)^+}e\le\frac{1}{\epsilon^2}T(f^2).$$
 \end{thm}

\proof  Let $f\in {\cal L}^2(T)$.
 As $P_{(f-\epsilon e)^+}$ is the band projection onto the band generated by
 $(f-\epsilon e)^+$ it follows that $P_{(f-\epsilon e)^+}(f-\epsilon e)\ge 0$ and thus
 \begin{eqnarray}
  P_{(f-\epsilon e)^+}f\ge \epsilon P_{(f-\epsilon e)^+}e\ge 0.\label{2014-1}
 \end{eqnarray}
 Band projections are dominated by the indentity map, so $P_{(f-\epsilon e)^+}\le I$, giving
 $|f|\ge P_{(f-\epsilon e)^+}|f|.$
From the positivity of band projections, $P_{(f-\epsilon e)^+}|f|\ge P_{(f-\epsilon e)^+}f$.
Taking these observations together with (\ref{2014-1}) gives
 \begin{eqnarray}
  |f|\ge \epsilon P_{(f-\epsilon e)^+}e\ge 0.\label{2014-2}
 \end{eqnarray}
 Multiplying (\ref{2014-2}) successively by $|f|, \epsilon P_{(f-\epsilon e)^+}e \ge 0$ in
 the $f$-algebra $E^u$, universal completion of $E$, gives
 \begin{eqnarray}
  f^2=|f|^2\ge \epsilon |f|P_{(f-\epsilon e)^+}e\ge (\epsilon P_{(f-\epsilon e)^+}e)^2\ge 0.\label{2014-3}
 \end{eqnarray}
 The construction of the $f$-algebra structure on $E^u$ yields directly that
 $Qe\cdot Qe=Qe$ for all band projections $Q$. Hence
 \begin{eqnarray}
   (P_{(f-\epsilon e)^+}e)^2=P_{(f-\epsilon e)^+}P_{(f-\epsilon e)^+}e^2=P_{(f-\epsilon e)^+}e.
    \label{2014-4}
 \end{eqnarray}
 Combining (\ref{2014-3}) and (\ref{2014-4}) gives
 \begin{eqnarray}
  f^2\ge \epsilon^2 P_{(f-\epsilon e)^+}e\ge 0.\label{2014-5}
 \end{eqnarray}
 Noting that $f^2\in E$,  $T$ can be applied to (\ref{2014-5}) to give the desired inequality.
\qed

%%%%%%%%%%%%%%%%%%%%%%%%%%%
\newsection{Bienaym\'e Equality}

The Bienaym\'e equality of classical statistics gives that the variance of a finite sum of 
independent random variables coincides with variance of their sum.  In this section we give
a measure free conditional version of this result in  ${\cal L}^2(T)$. 
Before we can proceed with this we require a result on $T$-conditionally independent random
 variables in ${\cal L}^2(T)$.

\begin{lem}\label{lemma}
Let $E$ be a $T$-universally complete Riesz space with weak order unit, $e=Te$, where $T$ is a 
strictly positive conditional expectation operator on $E$.
Let $f, g \in {\mathcal L}^2(T)$.  If $f$ and $g$ are $T$-conditionally independent then 
	\[Tfg = Tf\cdotp Tg = Tg \cdotp Tf.\]
\end{lem}

\proof
Let $T_f$ and $T_g$ denote the conditional expectations with ranges 
$\left< {\cal R}(T), f\right> = E_f$ and  $\left< {\cal R}(T), g\right> = E_g$ respectively. 
Here $\left< {\cal R}(T), g \right>$ denotes the order
 closed Riesz subspace of $E$ generated by ${\cal R}(T)$ and $g$, and similarly for
$\left< {\cal R}(T), f\right>$.
 The existence and uniqueness of $T_f$ and $T_g$ are given by the 
Radon-Nikod\'{y}m Theorem, see \cite{watson}. 
Here $E_g$ and $E_f$ are $T$-conditionally independent as the $T$-conditional independence of 
$f$ and $g$ is defined in terms of the independence of $E_f$ and $E_g$, see Section 2.
For each $h\in E$, $T_fh\in E_f$ and  $T_gh\in E_g$. 
Now,  as $E_f$ and $E_g$ are $T$-conditionally independent, 
from (\ref{eq-indep-1}) with $E_1=E_f$ and $E_2=E_g$, we have
\begin{equation}
	T_f(T_gh) = Th = T_g(T_fh). \label{ave prop}
\end{equation}
Applying (\ref{ave prop}) with $h=fg$ gives
\begin{align*}
T(fg) &= T_fT_g(fg)
\end{align*}
As $T_f$ and $T_g$ are averaging operators in ${\cal L}^2(T)$,  by Corollary \ref{cor-average},
$T_g(fg) = gT_gf$. Thus
\begin{align}
T(fg) &= T_fT_g(fg)=T_f(gT_gf),\label{baw-jan-1}
\end{align}
however taking (\ref{ave prop}) with $h=f$ yields $T_gf=Tf$, which along with
 (\ref{baw-jan-1}) gives 
\begin{align}
T(fg)&=T_f(gTf).\label{baw-jan-2}
\end{align}
In (\ref{baw-jan-2}), $Tf\in {\cal R}(T)\subset E_f$, so by Corollary \ref{cor-average},
$T_f(gTf)=Tf\cdot T_f g$. Finally considering (\ref{ave prop}) with $h=g$ 
gives $T_fg=Tg$. Thus
$$T(fg)=T_f(gTf)=Tf\cdot T_f g=Tf\cdot T g.$$
\qed

From Theorem \ref{lemmastar}, if $f\in {\cal L}^2(T)$ then $Tf \in {\mathcal L}^2(T)$ which gives $(f - Tf)  \in {\mathcal L}^2(T)$.  
Hence, $(f - Tf)^2  \in {\mathcal L}^1(T)$  and so $T(f - Tf)^2$ exists for all $f \in {\mathcal L}^2(T)$.  We now define the variance 
of $f$ by 
\begin{equation}
\var(f) = T(f- Tf)^2 = Tf^2 - (Tf)^2. \label{cross}
\end{equation}

\begin{thm}[Bienaym\'e Equality]\label{bienayme}
Let $E$ be a $T$-universally complete Riesz space with weak order unit, $e=Te$, where $T$ is a 
strictly positive conditional expectation operator on $E$.
If $(f_k)_{k\in\N}$, is a $T$-conditionally independent sequence in ${\mathcal L}^2(T)$, then 
	\[\var\left(\sum_{k=1}^n f_k\right) = \sum_{k=1}^n\var(f_k),\]
for each $n\in\N$.
\end{thm}

\proof
As 
	\[\left< f_{i_1},\dots f_{i_j}, {\cal R}(T)\right> = \left< f_{i_1}- Tf_{i_1},\dots f_{i_j}-Tf_{i_j}, {\cal R}(T)\right>, \]
for each subset $\{i_1,\dots,i_j\}$ of $\{1,\dots,n\}$.
Thus $f_k - Tf_k$, $k = 1,2, \dots, n$ are $T$-conditionally independent.
 From Theorem \ref{lemmastar}, as $f_i \in {\mathcal L}^2(T)$ it follows that
 $Tf_i \in {\mathcal L}^2(T)$ and consequently from Lemma \ref{lemma} that 
	\[T[(f_{i}-Tf_{i})(f_{j}-Tf_{j})]=[T(f_{i}- Tf_{i})]\cdotp[T(f_{j}-Tf_{j})],\]
 for $i\ne j$.
However, as $T$ is a projection, see Section 2,
	\[T(f_k - Tf_k) = 0,\quad\mbox{for each}\quad k\in\N,\]
giving
\begin{align}
T[(f_{i}- Tf_{i})(f_{j}-Tf_{j})]=0, \label{plus}
\end{align}
for $i\ne j$.
From the definition of variance
\begin{eqnarray*}
\var\left(\sum_{k=1}^nf_k\right) = T\left(\sum_{k=1}^nf_k - T\sum_{k=1}^n f_k\right)^2
 =T\left(\sum_{k=1}^n(f_k - Tf_k)\right)^2,
\end{eqnarray*}
which can be expanded to give
\begin{eqnarray}
\var\left(\sum_{k=1}^nf_k\right)
	= T\sum_{k=1}^n (f_k - Tf_k)^2 + T\sum_{j\neq k} (f_j - Tf_j)(f_k - Tf_k)\label{2014-a}
\end{eqnarray}
Now applying (\ref{plus}) to (\ref{2014-a}) gives 
\begin{eqnarray*}
\var\left(\sum_{k=1}^nf_k\right) =   \sum_{k=1}^n T(f_k - Tf_k)^2
	= \sum_{k=1}^n\var(f_k).
\end{eqnarray*}
 \qed

%%%%%%%%%%%%%%%%%%%%%%%%%%%
\newsection{Bernoulli and Poisson Processes}

 In classical probability, a Bernoulli process is one in which the events at any given time are 
 independent of the events at all other times.
 The payoff of an event occuring is $1$ unit and $0$ units for it not occuring.
 Thus in the Riesz space setting, the process can be described by the sequence of 
 independent band projections $P_k$ where $k$ indexes time and 
 the payoff at time $k$ is $P_ke$.  The probability of an event at time $k$ occuring 
 must be independent of $k$, in the measure theoretic terms, this can be expressed as 
 the expectation of each event is independent of time.
 This can be generalized to the conditional expectation of the events being time invariant,
 which lead to the Riesz space setting requirement that $TP_ke=f$, for all $k\in\N$.
 Here $T$ is some fixed conitional expectation operator. Thus we are led to the
 following formal definition of a Bernoulli process in Riesz spaces.
 
 \begin{defs}  
 Let $E$ be a Dedekind Riesz space with weak order unit, $e$, and conditional expectation operator 
 $T$ with $Te=e$.
 Let $({ P}_k)_{k\in\N}$ be a sequence of $T$-conditionally 
 independent band projections.  We say that $(P_k)_{k\in\N}$ is a Bernoulli
 process if 
 	\[ TP_ke = f \quad \mbox{for all }k \in \N,\]
 for some fixed $f\in E$.
\end{defs}

The payoff at time $n$ is thus 
 	\[S_n = \sum_{j=1}^n P_je.\]
We denote by $P_{S_n = je}$ the band projection on the band where $S_n = je$, in the notation
used earlier $P_{S_n = je}=(I-P_{(S_n - je)^+})(I-P_{(S_n - je)^-})$.

\begin{thm}\label{bernoulli}
Let $E$ be a $T$-universally complete Riesz space with weak order unit, $e=Te$, where $T$ is a 
strictly positive conditional expectation operator on $E$.
Let $(P_j)_{j\in\N}$ be $T$-conditionally independent band projections with $TP_je =f$ 
for all $j\in\N$ and 
$S_n =\displaystyle\sum_{j=1}^nP_je$.  Then
\begin{eqnarray}
        TS_n&=&nf,\label{mean}\\
	TP_{S_n = je}e&=&\frac{n!}{j!(n-j)!}f^j(e-f)^{n-j},\label{expectation}\\ 
	\var(S_n)&=&nf(e-f).\label{variance}
\end{eqnarray}
\end{thm}

\proof
 As $TP_{i}e = f$, (\ref{mean}) follows directly from applying $T$ to $S_n$.

 Fix $n\in\N$ and 
 let $$Q_j = \frac{1}{j!(n-j)!}\sum_{\sigma\in\Lambda}
   P_{k_{\sigma (1)}}\dots P_{k_{\sigma (j)}}(I-P_{k_{\sigma (j+1)}})\dots (I-P_{k_{\sigma (n)}}).
   \quad j=0,\dots,n.$$ 
 Here $\Lambda$ denotes the set of all permutations of $\{1,\dots,n\}$ and the division by
  $j!(n-j)!$ is as there are $j!(n-j)!$ permutations which yield the same band projection
 $P_{k_{\sigma (1)}}\dots P_{k_{\sigma (j)}}(I-P_{k_{\sigma (j+1)}})\dots (I-P_{k_{\sigma (n)}})$.
 Other permutations yield band projections disjoint from the above one.  Thus $Q_j$ is a 
 band projection, $Q_0,\dots, Q_n$ partition the identity, $I$, in the sense that
  $Q_iQ_j=0$ for all $i\ne j$, and $\sum_{i=0}^n Q_i=I$.
 Moreover, from the definition of $Q_j$, it follows that 
 $Q_jS_n=jQ_je, j=0,\dots,n.$  Thus
 	\[S_n=\sum_{j=0}^n Q_jS_n = \sum_{j=0}^n j Q_je.\]
The $T$-conditional independence of $P_1, \dots, P_n$ and Lemma \ref{lemma} applied iteratively give that 
\begin{align*}
&TP_{k_{\sigma (1)}}\dots P_{k_{\sigma (j)}}(I-P_{k_{\sigma (j+1)}})\dots (I-P_{k_{\sigma (n)}})e\\
&=T((P_{k_{\sigma (1)}}\dots P_{k_{\sigma (j)}}e)\cdot
       ((I-P_{k_{\sigma (j+1)}})\dots (I-P_{k_{\sigma (n)}})e))\\
&=(T(P_{k_{\sigma (1)}}\dots P_{k_{\sigma (j)}}e)\cdot
       (T(I-P_{k_{\sigma (j+1)}})\dots (I-P_{k_{\sigma (n)}})e))\\
&  =\prod_{i=1}^jTP_{k\sigma(i)}e\cdot\prod_{i=J+1}^nT(I-P_{k_{\sigma (i)}})e\\
&  =f^j(e-f)^{n-j}.
\end{align*}
Hence
	\[TP_{S_n = je}e = TQ_je = \frac{1}{j!(n-j)!}\sum_{\sigma\in\Lambda}f^j(e-f)^{n-j},\]
from which (\ref{expectation}) follows as the cardinality of $\Lambda$ is $n!$.

As $P_1e,\dots, P_ne$ are $T$-conditionally independent and are in $\mathcal{L}^2(T)$, 
Bienaym\'e\rq{}s equality applied to $S_n$ gives
\begin{eqnarray}
  \var(S_n) = \sum_{k=1}^n \var(P_ke).\label{eq-var}
\end{eqnarray}
From (\ref{cross}) applied to $P_ke$ we have
\begin{eqnarray}
\var(P_ke) =TP_ke - (TP_ke)^2
	= f -f^2.\label{var-2}
\end{eqnarray}
As $e$ is the multiplicative unit, combining (\ref{eq-var}) and  (\ref{var-2}) yields (\ref{variance}).
\qed

\begin{thm}[Bernoulli Law of Large Numbers]
Let $E$ be a $T$-universally complete Riesz space with weak order unit, $e=Te$, where $T$ is a 
strictly positive conditional expectation operator on $E$.
Let $(P_k)_{k\in\N}$ be a Bernoulli process with partial sums $S_n$ and
$TP_ke=f, k\in\N$.
For each $\epsilon > 0$,  
	\[TP_{\left(\vert\frac{S_n}{n} - f\vert -\epsilon e\right)^+}e \to 0,\]
as $n \to\infty$.
\end{thm}

\proof
By the Tchebichev inequality, 
\begin{eqnarray}
TP_{\left(\vert S - nf\vert -n\epsilon e\right)^+}e 
	\le \frac{1}{n^2\epsilon^2}T\vert S_n-nf\vert^2.\label{blln-1}
\end{eqnarray}
However, from (\ref{cross}) and (\ref{mean}),
\begin{eqnarray}
	T\vert S_n-nf\vert^2=\var(S_n).\label{blln-2}
\end{eqnarray}
Combining (\ref{blln-1}) with (\ref{blln-2}) and using (\ref{variance}) to simplify the result, gives
\begin{eqnarray}
TP_{\left(\vert S - nf\vert -n\epsilon e\right)^+}e 
	\le \frac{f(e-f)}{n\epsilon^2},
\end{eqnarray}
from which the result follows upon observing that
 $P_{\left(\vert\frac{S_n}{n} - f\vert -\epsilon e\right)^+}=P_{\left(\vert S - nf\vert -n\epsilon e\right)^+}$.
\qed

One of the interesting features of Bernoulli's law of large numbers is that it gives not just the 
convergence of $TP_{\left(\vert\frac{S_n}{n} - f\vert -\epsilon e\right)^+}e$ to zero. 
It also gives some 
indication of the size of the band on which $\vert\frac{S_n}{n} - f\vert > \epsilon e$
by bounding the conditional expectation of the band projection applied to $e$ 
 by $\frac{f(e-f)}{n\epsilon^2}$, hereby indicating an upper bound for the rate of convergence
`in probability'.

Using the results on martingale difference sequences developed for the study of mixingale in
\cite[Lemma 4.1]{kvw} we obtain a weak law of large numbers for Bernoulli processes.

\begin{thm}[Weak law of large numbers]
Let $E$ be a $T$-universally complete Riesz space with weak order unit, $e=Te$, where $T$ is a 
strictly positive conditional expectation operator on $E$.
Let $(P_j)_{j=1,\dots, n}$ be $T$-conditionally independent band projections with $TP_je =f$ 
for all $j =1, \dots, n$ and 
$S_n =\displaystyle\sum_{j=1}^nP_je$, then
 \begin{eqnarray*}
     \lim_{n\to\infty} T\left|f-\frac{S_n}{n} \right|=0.
 \end{eqnarray*}
\end{thm}

\proof
 Setting $f_i=P_ie$ and $T_i$ to be the conditional expectation with range 
 $\left< {\cal R}(T), P_1e, \dots, P_ie\right>$, it follows that $(g_i,T_i)$, where
 $g_i:=f_i-T_{i-1}f_i$, is a martingale difference sequence.
 Here $|f_i|\le e$. 
 Thus from \cite[Lemma 4.1]{kvw},
 \begin{eqnarray}
     \lim_{n\to\infty} T\left|\frac{1}{n}\sum_{i=1}^n  g_i   \right|=0.
    \label{wll}
 \end{eqnarray}
 The independence of the band projections $P_i, i\in\N,$ gives that 
 $T_{i-1}f_i=Tf_i=f$.  Hence (\ref{wll}) can be written as
 \begin{eqnarray*}
     \lim_{n\to\infty} T\left|f-\frac{1}{n}\sum_{i=1}^n  P_ie   \right|=0,
 \end{eqnarray*}
 from which the theorem follows.
\qed

Before progressing further we need to define an 
exponential map on Riesz spaces.

\begin{rem}\label{rem} {\rm
Let $C([-1,1])$ denote the Riesz space of continuous real functions on $[-1, 1]$.
Set $f_n(t) := \displaystyle\left(1-\frac{t}{n}\right)^n$, then $f_n \in C([-1,1])$ 
and $f_n(t) \to e^{-t} = f(t)$ in order and 
 the supremum norm on $C([-1,1])$.  
Thus, by \cite[Theorem 3.1]{grobler2}, for each $g \in E^e, f_n(g) \to f(g)$ $e$-uniformly 
 (and, thus, in order) as $n \to \infty$.  In addition, by the functional calculus, $f(g)$ defines an element of $E^e$ which we will
 denote by $e^{-g}$. }
\end{rem}

We now consider the sequences of 
Bernoulli processes known as Poisson sequences. Here the partial sums of each Bernoulli process
form a Bernoulli process.

\begin{thm}[Poisson]\label{poisson}
Let $E$ be a $T$-universally complete Riesz space with weak order unit, $e=Te$, where $T$ is a 
strictly positive conditional expectation operator on $E$.
Let $P_{n,k}$, $k=1,\dots, n,\  n \in \N$, be $T$-conditionally independent band projections
 with $TP_{n,k}e = g_n$ for all $k =1, \dots, n, \ n\in \N$.
If $\displaystyle S_n =\sum_{k=1}^nP_{n,k}e$ are $T$-conditionally independent
with $TS_n = g$, $n\in \N$, then for each $j=0,1,\dots$,
	\[TP_{S_n = je}e\to \frac{g^j}{j!}e^{-g},\]
$e$-uniformly as $n\to\infty$.
\end{thm}

\proof
 From (\ref{expectation}), 
 \begin{eqnarray*}
   TP_{S_n = je}e =\frac{n!}{j!(n-j)!}g_n^j(e-g_n)^{n-j},
 \end{eqnarray*}
 but (\ref{mean}) gives $g=TS_n = ng_n$.  Hence 
 \begin{eqnarray*}
   TP_{S_n = je}e =\frac{n!}{j!(n-j)!}\left(\frac{g}{n}\right)^j\left(e-\frac{g}{n}\right)^{n-j},
 \end{eqnarray*}
 which can be expanded to give
 \begin{eqnarray}
   \left(e-\frac{g}{n}\right)^{j}TP_{S_n = je}e =\frac{g^j}{j!} 
\left(1-\frac1n\right)\left(1-\frac2n\right)\cdots\left(1-\frac{j-1}n\right)
\left(e-\frac{g}{n}\right)^{n}.\label{poisson-5}
 \end{eqnarray}
 Taking the limit as $n\to \infty$ in (\ref{poisson-5}) gives
\begin{align*}
 TP_{S_n = je}e= \frac{g^j}{j!}\lim_{n\to\infty}\left(e -\frac{g}{n}\right)^n,
\end{align*}
 which together Remark \ref{rem} concludes the proof. 
\qed

%%%%%%%%%%%%%%%%%%%%%%% bibliography %%%%%%%%%%%%%%%%%%

\end{document}